\documentclass[11pt]{article}
\usepackage[letterpaper]{geometry}
\usepackage{amsthm,amsmath,amsfonts}
\usepackage[normalem]{ulem}
\usepackage{color,graphicx}
\usepackage{float}
\usepackage[colorlinks]{hyperref}

\usepackage[numbers]{natbib}
\usepackage{centernot}
\usepackage{enumerate}
\newcommand{\R}{\mathbb{R}}

\newcommand{\E}{\mathbb{E}}

\newcommand{\A}{\mathcal{A}}

\newcommand{\Prob}{\mathbb{P}}
\usepackage{graphicx}
\graphicspath{{images/}{../images/}}

\newtheorem{lemma}{Lemma}

\theoremstyle{definition}
\numberwithin{equation}{section}

 \usepackage{thmtools, thm-restate}

\usepackage{tikz}
 \usepackage{tkz-graph}
 \usetikzlibrary{arrows}
\usetikzlibrary{positioning}
\usetikzlibrary{calc,fit}
\usepackage[font=small, labelfont=bf]{caption}
\usepackage{subfig}
\usepackage{float}
\usepackage[ruled,vlined,linesnumbered]{algorithm2e}
\usepackage{xcolor}
\usepackage{authblk}

\usepackage[utf8]{inputenc}
\usepackage{quoting}
\usepackage{graphicx}
\usepackage{tikz}
\usetikzlibrary{positioning}
\usetikzlibrary{calc,fit}
\usepgflibrary{shapes.geometric}
 
\begin{document}

\title{Scalable Deep Reinforcement Learning for Ride-Hailing}
\date{\vspace{-5ex}}

\author[1]{Jiekun Feng}
\author[2]{Mark Gluzman}
\author[,3,4]{J. G. Dai\footnote{Corresponding author:
J. G. Dai,
email: jd694@cornell.edu}}
\affil[1]{Department of Statistics and Data Science, Cornell University, Ithaca, NY}
\affil[2]{Center for Applied Mathematics, Cornell University, Ithaca, NY}
\affil[3]{ School of Data Science,  Shenzhen Research Institute of Big Data, The Chinese
University of Hong Kong, Shenzhen}
\affil[4]{School of Operations Research and Information Engineering, Cornell University, Ithaca, NY}

\maketitle
 \abstract{
 Ride-hailing services, such as Didi Chuxing, Lyft, and Uber,  arrange thousands of cars to meet ride requests throughout the day. We consider a Markov decision process (MDP) model of a ride-hailing service system, framing it as a reinforcement learning (RL) problem. The simultaneous control of  many agents (cars) presents a challenge for the MDP optimization because the action space grows exponentially with the number of cars. We propose a special decomposition for the MDP actions by sequentially assigning tasks to the drivers. The new actions structure resolves the scalability   problem and enables the use of deep RL algorithms for  control policy optimization. We demonstrate the benefit of our proposed   decomposition with a numerical experiment based on real  data from Didi Chuxing.
 }

\section{Introduction} 
\label{sec:intro}

A ride-hailing service  is the next generation of taxi service that uses online platforms and mobile apps to connect passengers and drivers. Lyft, Uber, and Didi Chuxing together serve more than 45 million passengers per day \cite{Schlobach2018}.

 A centralized planner of the ride-hailing service arranges cars in the system, matching them with ride requests from the passengers.  The centralized planner may also relocate an empty (without a passenger) car to another location in anticipation of future demand and/or shortage.  In other words, the centralized planner   assigns tasks to multiple drivers over a certain time horizon, considering the expected geographical distributions of  requests and cars.

  Motivated by the problem of relocating cars with no passengers, we follow \cite{BravDaiLiuLei2019}, which used a closed queueing network model under the assumption of constant traffic   parameters (e.g. passengers arrival rates, cars travel times),
 formulated a fluid-based optimization
problem and found an asymptotically optimal empty-car
routing policy in the “large market” parameter region. For
a time-varying traffic pattern, each decision time the traffic parameters were averaged over the finite time window and used to formulate the fluid-based
optimization problem  as if the ride-hailing service had a constant traffic pattern.  The result was a time-dependent
lookahead policy that  was not designed to be optimal.

  The reinforcement learning (RL)   framework we propose, however, assumes that the centralized planner does not know the traffic parameters. An RL problem refers to a (model-free) Markov decision process (MDP) problem in which the underlying dynamics is unknown, but optimal actions can be learnt from sequences of data (states, actions, and rewards) observed or generated under a given policy. 
  The centralized planner  has to  control multiple agents (cars) simultaneously while accounting for their interactions. Since the centralized planner's action space will grow exponentially with the number of agents, we propose a method that transforms any action of the centralized
planner into a special form that enables both simultaneous control of multiple agents and the use of policy
optimization RL algorithms \cite{Kakade2002, Schulman2015, Schulman2017} to optimize the centralized planner's decisions.

The development and application of RL algorithms combined with neural networks for controlling and optimizing ride-hailing services are important research topics. 
Related to empty-car routing,  \cite{Oda2018} employed a   deep Q-network (DQN)  algorithm, proposed in \cite{Mnih2015}, to proactively dispatch cars to meet  future  demand while limiting   empty-car routing times. To ensure scalability, the algorithm was applied to find the optimal dispatch actions for individual cars, but without considering agents' interactions. We note that the use of 
deep RL algorithms to learn the optimal actions from the perspective of individual cars has been studied intensively \cite{Oda2018, WangQinTanYeZhu18, ShiGao19}, yet this approach can be ``intrinsically inadequate to be used in a production dispatching system which demands coordinations among multiple agents'' \cite{Tang2019ADV}.
In~\cite{Tang2019ADV} the authors studied ``ride-sharing order dispatching'' processes using a two-stage optimization procedure consisting of  a deep RL algorithm followed by combinatorial optimization for bipartite graph matching. A simulation study based on real dispatch data from Didi Chuxing showed that the authors' deep RL approach outperformed other recently proposed dispatching methods ~\cite{WangQinTanYeZhu18, Xu2018ADV}   in terms of driver income and passenger experience metrics.
In~\cite{KeXiao19} the authors proposed to delay the car assignment for some passengers in order to   accumulate more   drivers and waiting passengers in the matching pool, where the passengers and cars would be connected later via optimal bipartite matching.  
A multi-agent
actor-critic RL algorithm was developed to optimize the choice of  the delayed times for the passengers. 
The use of Didi Chuxing's real data demonstrated that the proposed actor-critic method decreased the mean pickup time significantly compared to pure optimization methods that assumed no delay.

In this paper we consider a centralized planner that receives real-time data on existing ride
requests and driver activities and assigns tasks to drivers at each decision epoch. We assume the
decision epochs occur at discrete times, and the time between every two consecutive decision epochs is fixed. At
each decision epoch the centralized planner must solve the following combinatorial problem:   each available car should be either  matched with a passenger’s ride request, relocated to another location with no passengers, or asked to stay at its current location until the next decision epoch.

This paper makes the following contributions to the ride-hailing literature. The new method proposed to solve combinatorial problems uses a special decomposition of actions  and a control policy that \emph{sequentially} matches agents with tasks. A deep learning algorithm is used to optimize the control policy  that should generate the most beneficial trip to fulfill at a current system's state. The centralized planner's control policy that assigns tasks is based on mapping from the state space to a set of trips.  The complexity of its representation grows linearly with the ride-hailing service area. We use a proximal policy optimization (PPO) algorithm \cite{Schulman2017} for the control policy optimization. The results of numerical experiment using real data from a ride-hailing service show that the proposed solution outperforms currently used combinatorial matching method.

  The remainder of this paper is organized as follows. Section \ref{sec:ridesharing-network} describes the ride-hailing service and its transportation network. Section \ref{sec:prob-form} formulates the optimal control problem and control optimization algorithm. Section \ref{sec:num-study} compares the results of the proposed method to the  time-dependent lookahead policy proposed by  \cite{BravDaiLiuLei2019} based on a numerical study using real data from a ride-hailing service. Section \ref{sec:conclusion} concludes and offers suggestions for future research.

\section{The transportation network}
\label{sec:ridesharing-network}
 
 In this section we describe our model of the ride-hailing service and transportation network, following \cite{BravDaiLiuLei2019}, with a few modifications. The service consists of a centralized planner, passengers requesting rides, and a fixed number of geographically distributed agents (cars). The transportation network consists of $N$ cars distributed across a service territory divided into $R$ regions.
For   ease of exposition, we assume that each working day (``episode'') of the  ride-hailing service starts at the same time  and  lasts for $H$ minutes.

  We assume that the number of passenger arrivals at region $o$ in the $t$-th minute (i.e., $ t$ minutes elapsed since the start of the working day) is a Poisson random variable with mean
  \begin{align*}
    \lambda_o(t),  \text{ for each } o=1,\dotsc,R, \, t=1,\dotsc,H.
\end{align*}
 The collection of all  Poisson random variables is independent. Passengers only arrive \emph{after} a working day starts (i.e., there are no passengers at the 0-th minute).
  


 
Upon arrival  at region $o$,  a passenger travels to region $d$ with probability that depends on  time $t$, origin region $o$, and destination region $d$
\begin{align*}
    P_{od}(t), \quad o,d=1,\dotsc,R, \, t=1,\dotsc,H.
\end{align*}
 
After a trip   from region $o$ to $d$ has been initiated, its duration is \emph{deterministic} and equals to 
\begin{align}
    \tau_{od}(t), \quad o,d=1,\dotsc,R, \, t=1,\dotsc,H.
    \label{eq:trav-time}
\end{align}

  We let 
 \begin{align}\label{eq:tau_d}
 \tau_d :=\max_{t=1,\dotsc,H, \, o=1,\dotsc,R}\tau_{o d}(t),\quad   d=1,..., R
 \end{align} be the maximum travel time to  region $d$ from any region of the transportation network at any time.

 While Section 1  in \cite{BravDaiLiuLei2019}  assumed that travel times were random variables having
an exponential distribution, the experiments in Section 3.2.1 were conducted under constant travel times. For ease of exposition, we use deterministic travel times in Section \ref{sec:prob-form} below.


Patience time denotes a new passenger's maximum waiting time for a car. We assume that each passenger  has a deterministic patience time   and we fix it as equal to $L$ minutes.  We assume that the centralized planner knows the patience time.  
  
In real time, the centralized planner receives   ride requests, observes the location and activity of each car in the system, and considers three types of tasks for the available cars: (1) car-passenger matching, (2) empty-car routing, and  (3)  ``do nothing''  (a special type of   empty-car routing).  
  We assume that each passenger requires an immediate response to a request within the first decision epoch.   If the centralized planner assigns a matching between a passenger and an available car, we assume  the passenger has to accept the matching. A passenger who is not matched with a car in the first decision epoch leaves the system.

Unlike \cite{BravDaiLiuLei2019}, we relax the constraint that   only  cars   idling at the passenger's location can be matched with the passenger. We assume that the centralized planner can match cars with subsequent ride requests before current trips are completed.
We assume that the   patience time   satisfies 
 \begin{align} 
    L&<\min_{t=1,\dotsc,H, \, o=1,\dotsc,R}
    \tau_{o d}(t), 
    \quad \text{for each } d = 1,\dotsc,R
    .
    \label{eq:a-trav-time}
\end{align}
The assumption implies that the travel time of any trip   is larger than the patience time.    Therefore, no more than one subsequent trip can be assigned to a driver.

If a car  reaches its destination and has not been matched with a new passenger, it becomes empty.
The centralized planner may let the empty car stay at the destination or relocate to another region.  For the former, we note that the centralized planner's decision belongs to  the ``do nothing'' task and does not cost any travel time. The centralized planner will be able to assign the car a new task at the next decision epoch. For the latter, the centralized planner chooses a region  for the relocation and the travel time remains the same as in  Equation~\eqref{eq:trav-time}.   
 Unlike \cite{BravDaiLiuLei2019}, the centralized planner can assign two empty-car routing tasks in succession.


\section{Optimal control problem formulation}
\label{sec:prob-form}

 Our goal is to find a control policy for the centralized planner that maximizes the total reward collected during one working day by the entire ride-hailing service. In this section we formulate the problem as a   finite-horizon, discrete-time,  undiscounted MDP. We set the time interval between two successive epochs to one minute.
As a result, a passenger waits at most one minute for a decision. Under this setting, the time in minute, $t = 1,\dotsc,H$, also represents the decision epochs. 
 
\subsection{State space}
\label{subsec:state}

 The state space $S^\Sigma$ of the MDP includes states $s_t = \left[ s_t^e, s_t^c, s_t^p\right]$, such that each state consists of three components:   current epoch  $s_t^e := t$,  cars status  $s_t^c$, and   passengers status  $s_t^p$. 
 
The cars status component represents the  number of   cars of every   type   in the system:
\begin{align*}
s_t^c &:= \left(s_t^c(d,\eta )~\Big| ~  d = 1,\dotsc, R, ~ \eta = 0,1,\dotsc,\tau_d, \tau_d+1,\dotsc,\tau_d+L \right), 
\end{align*}
 where $s_t^c(d,\eta )$ is the number of cars in the system whose \emph{final} destination region is $d$,  and the \emph{total} remaining travel time  (``distance'') to the destination is equal to $\eta$, and $\tau_d$ is the maximum travel time to region $d$ defined by (\ref{eq:tau_d}).
 
The passengers status component is equal to
\begin{align*}
s^p_t &:= \left( s^p_t(o, d)~\Big|~o, d=1,..., R \right),
\end{align*}
where $s^p_t(o, d)$ characterizes the number of passengers in the system requesting rides from region $o$ to region $d$, $o, d = 1, ..., R.$

\subsection{Action space  }
\label{subsec:action}

At each epoch $t$, the centralized planner  observes  the system state $s_t$, and makes a decision $a_t$   that should  address  all    $I_t$ available cars, where
\begin{align*}
   I_t:= \sum_{o = 1}^R \sum_{\eta = 0}^L s_t^c (o, \eta).
\end{align*}

 We let $\A^\Sigma$ denote the action space of the MDP. 
  We propose to  decompose every decision $a_t\in \A^\Sigma$ into a sequence of ``atomic actions'', each addressing a single available car, to overcome the challenge of the large action space. Therefore, we consider action $a_t$ as:
\begin{align*}
  a_t := \left(  a_{t,1}, ..., a_{t, I_t}\right),
\end{align*}
where $a_{t, i}$ is an atomic action that encodes a trip   by one of the candidate cars. We let $\A$  denote the atomic action space. We  note that $\A = \{(o,d)\}_{o,d=1}^R$.

We call the  sequential generation of atomic actions a ``sequential decision making process'' (SDM process). We let $s_{t,i}$   denote  a state of the SDM process after   $i-1$ steps, for each decision epoch $t=1,..., H$. We let $S$ be the state space of the SDM process.
Each  state $s_{t,i}$ of the SDM process has four components $s_{t,i} := [s_{t,i}^e, s_{t, i}^c, s_{t,i}^p, s_{t, i}^\ell]$, where, as in the original MDP, the first three components $s_{t,i}^e$, $s_{t,i}^c$, $s_{t, i}^p$  represent current epoch,  cars status, and passengers status, respectively, and a new component $s_{t,i}^\ell$ tracks the cars exiting the available cars pool until the next decision epoch.
The SDM process is initialized with state $s_{t, 1}$ such that 
$s_{t,1}^e = s_t^e,$ $s_{t, 1}^c = s_t^c,$ $s_{t,1}^p= s_t^p$,
 and   $s_{t,1}^\ell$ is  a zero vector,  for each decision epoch $t=1,...,H$.

Each atomic action represents a ``feasible'' trip $a_{t,i}  = (o_{t,i}, d_{t,i})$, 
where $o_{t,i},$ $d_{t,i}$ are the origin   and  destination regions of the trip, respectively.
Action $a_{t,i}$  is \emph{feasible} if there exists an available car that is $L$ minutes (or less) away from the origin region $o_{t,i}$, (i.e. $\sum_{\eta = 0}^{L} s^c_{t,i}(o_{t,i}, \eta) > 0$).
 Although atomic action $a_{t,i}$  only encodes the   origin and destination of a trip, we set a few rules that   specify which car will conduct the trip and if the car will carry a passenger. 
 Among the set of available cars, we select  the car closest to origin $o_{t,i}$ to take the trip. 
Then,  we  prioritize car-passenger matching over empty-car routing, (i.e. if there exists a passenger requesting a ride from the trip origin to the trip destination, we assign the car to the requesting passenger; if there are several passengers requesting such a ride, we assign the car to a passenger at random).
If there is no passenger requesting a ride from $o_{t,i}$ to $d_{t,i}$, we interpret atomic action $a_{t,i}$ as either an empty-car routing task or a ``do nothing'' task.  If the available car assigned to the trip idles at origin region $o_{t,i}$ \emph{and} the trip relocates the car to a different region ($d_{t,i} \neq o_{t,i}$), then the car fulfills an empty-car routing task. Otherwise, we interpret the atomic action as a "do nothing" task, and the car becomes a "do-nothing" car. 

Once an available car,  possibly  a ``do nothing'' car,  has been assigned  a task at the $i$th step of the SDM process, the centralized planner should exclude it from the available cars pool. If the car has been assigned  a passenger ride request  or an empty-car routing task,   the cars status component of the SDM process state is updated such that the car becomes associated with its new final destination.
The car is automatically excluded from the available cars pool by assumption (\ref{eq:a-trav-time}). The ``do nothing'' tasks require   special transitions that the original MDP does not have, so we use $s_{t,i}^\ell$ to track the ``do nothing'' cars
\begin{align*}
s_{t,i}^\ell &= \Big(s_{t,i}^\ell(d,\eta) ~\Big|~ d = 1,\dotsc, R, ~\eta = 0,1,\dotsc,L \Big),
\end{align*}
 where $s_{t,i}^\ell(d,\eta)$ is the number of  ``do nothing'' cars which drive to or idle at region $d$, $\eta$ minutes away from their destinations.  
  The "do nothing" component excludes "do nothing" cars from the  available cars pool until the next decision epoch.




 The atomic actions are generated sequentially  under   control policy 
 \begin{align*}
 \pi:S\rightarrow \A,
 \end{align*}
 which is a mapping from the state space into a set of the trips. The  control policy $\pi$, given a current state of the SDM process $s_{t,i}$, sequentially generates feasible atomic actions $ a_{t,i} = \pi(s_{t,i})$. The   SDM process terminates when all candidate cars become non-candidate cars, producing action $a_t = (a_{t,1}, ..., a_{t,I_t})$.

 At each decision epoch the centralized planner observes system state $s_t$ and  exercises control policy $\pi$ sequentially  in  the SDM process to obtain action $a_t$. Then the transition of the system to the next state $s_{t+1}$ occurs according to the dynamics of the original MDP.

\subsection{Reward functions and objective}

A car-passenger matching generates an immediate reward
\begin{align*}
    c_{t}^f(o, d, \eta), \quad  o, d = 1,\dotsc,R, \,
    \eta = 0, 1, \dotsc, L,  ~
    t = 1,\dotsc, H,  
\end{align*}
where  $o$ and $d$ are the passenger's origin and destination regions, respectively, $\eta$ is the distance (in minutes) between the  matched car  and the passenger's location, and $t$ is the time of the decision.  The superscript ``f'' denotes a full-car trip.

Every  empty-car routing atomic action  generates  a cost that  depends on   origin region $o$, destination region $d$,   and   decision time $t$
\begin{align*}
     c_{t}^e(o,d), \quad o,~d = 1,\dotsc,R, \, t = 1,\dotsc, H, 
\end{align*}
where the superscript denotes an empty-car trip.

We assume the ``do nothing'' actions do not generate any rewards.
Therefore, a one-step reward function generated on  the $i-1$th step of SDM process at epoch $t$ is equal to 
\begin{align*}
     c(s_{t,i}, a_{t,i}) =  
     \begin{cases}
    c_t^f(o,d,\eta), \text{ if action }a_{t,i} \text{ implies car-passenger matching,}\\
    -c^e_t(o,d), \text{ if action }a_{t,i} \text{ implies empty-car routing},\\
    0, \text{ if action }a_{t,i} \text{ implies ``do nothing'' action}.
    \end{cases}
\end{align*}

We want to find   control policy $\pi$  that maximizes the expected total  rewards over the finite time horizon
\begin{align*}
    \E_\pi \left[\sum_{t = 1}^H \sum_{i = 1}^{I_t}
    c(s_{t,i}, a_{t, i})
    \right].
\end{align*}

\subsection{Control policy optimization}

Here,  a ``randomized'' control policy    refers to  a map $\pi:S \rightarrow [0,1]^{R^2}$, that given state $s\in S$ outputs a probability distribution over all trips. We use $\pi(a | s)$ to denote   the probability of choosing atomic action $a$ at state $s$ if the system operates under policy $\pi$. We assume that $\pi(a | s) =0 $ if action $a$ is infeasible at system state $s$. Then,  at each epoch $t$, the SDM process,  under randomized control policy $\pi$, selects atomic actions sampled according to   distribution $\pi(\cdot|s_{t,i})$, for each step $i=1,..., I_t.$

We define  a  value function $V_\pi: S\rightarrow \R$ of policy $\pi$  
\begin{align*} 
    V_\pi(s_{t,i}): =  \E_\pi   \left[\sum\limits_{k=i}^{I_t} c ( s_{t,k}, a_{t,k}) + \sum\limits_{j=t+1}^H\sum\limits_{k=1}^{I_j} c( s_{j,k}, a_{j,k}) \right],
\end{align*}
 for each $t=1, ..., H$,  $i=1,..., I_t$, and $s_{t,i}\in S.$  
 For notational convenience we set $V_\pi(s_{H+1, 1})=0$ for any policy $\pi$.

We define  advantage function $A_\pi:  S\times \A\rightarrow \R$ of policy $\pi$ 
\begin{align*} 
    A_\pi(s_{t,i}, a_{t,i}):=  
    \begin{cases}
    c(s_{t,i}, a_{t,i})+  V_\pi(s_{t, i+1}) - V_\pi( s_{t,i}), \quad \text{ if }i\neq I_t\\
     c(s_{t,i}, a_{t,i})+\sum\limits_{y\in S} \mathcal{P}(s_{t,i}, a_{t,i}, y) V_\pi(y) - V_\pi(s_{t, i}), \text{ if }i= I_t,\nonumber
    \end{cases}
\end{align*}
 for each $t=1,...,H$; $i=1,..., I_t$. We note that the transitions within the SDM process are deterministic. 
 We use $\mathcal{P}$ to denote the probabilities of transitions that come  from random  passenger arrivals.   We also define an advantage function for the original MDP $A^\Sigma_\pi: S^\Sigma\rightarrow \A^\Sigma$ of policy $\pi$  
 \begin{align*}
     A^{\Sigma}_\pi(s_{t}, a_t): = \sum\limits_{i=1}^{I_t} A_\pi(s_{t,i}, a_{t,i}),
 \end{align*}
  where $a_t = (a_{t,1}, ..., a_{t,I_t})$.  We note that $A_\pi^\Sigma (s_{t}. a_t) = \sum\limits_{i=1}^{I_t} c(s_{t,i}, a_{t,i}) + \sum\limits_{y\in S} \mathcal{P}(s_{t,I_t}, a_{t,I_t}, y)V_\pi(y) - V_\pi(s_{t,1})$. We also let $\pi^\Sigma(a_t|s_{t,1})$  denote  the  probability of selecting action $a_t = (a_{t,1}, ..., a_{t, I_t})$ under policy $\pi$ if the SDM process is initialized at state $s_{t,1}$.

 We let  $\{\pi_\theta,~ \theta\in \Theta\}$ be a set of parametrized control policies,
where $\Theta$ is an open subset of  $\R^q$,   $q\geq 1$.  
Hereafter, we abuse the notation and use $V_\theta$, $A_\theta$  to denote the value function and the advantage function  of policy $\pi_\theta$, $\theta\in \Theta$, respectively.

 In Lemma \ref{lem:per_diff} we obtain  a  performance difference equality for the MDP operating under the actions generated by the SDM process.    Similar performance difference equality was first obtained for MDPs with infinite-horizon discounted cost objectives in \cite{Kakade2002}.

 \begin{lemma}\label{lem:per_diff}
  We consider two policies $\pi_\theta$ and $\pi_\xi$, $\theta, \xi\in \Theta$.  
  Their value functions satisfy
 \begin{align*}
 V_\theta( s_{1, 1}) - V_\xi(s_{1, 1})  = \E_{\pi_\theta} \left[\sum\limits_{t=1}^{H} \sum\limits_{i=1}^{I_t} A_\xi ( s_{t,i}, a_{t,i})  \right].
  \end{align*}
 \begin{proof}
 
 First, we note that 
 \begin{align*}
 \E_{\pi_\theta} \left[V_\xi(  s_{1,1}) -V _\xi(  s_{H+1, 1}) + \sum\limits_{t=1}^{H }   \left( \sum\limits_{y\in S}\mathcal{P}( s_{t,I_t},a_{t,I_t}, y) V_\xi(y) - V _\xi(  s_{t,I_t}) + 
   \sum\limits_{i=1}^{I_t-1} \Big(    V_\xi(s_{t, i+1}) - V _\xi(  s_{t,i})\Big) \right)      \right]=0,
 \end{align*}
 see, for example, \cite{Henderson2002}.
 
 Then
\begin{align*}
&V_\theta(  s_{1, 1}) - V_\xi (s_{1, 1})  =  V_\theta( s_{1,1 }) - V_\xi( s_{1, 1}) \\
&\quad+ \E_{\pi_\theta} \left[V_\xi(  s_{1,1}) -V _\xi(  s_{H+1, 1}) + \sum\limits_{t=1}^{H }   \left( \sum\limits_{y\in S}\mathcal{P}( s_{t,I_t},a_{t,I_t}, y) V_\xi(y) - V _\xi(  s_{t,I_t}) + 
   \sum\limits_{i=1}^{I_t-1} \Big(    V_\xi(s_{t, i+1}) - V _\xi(  s_{t,i})\Big) \right)      \right]\\
     &= V_\theta(  s_{1, 1})  + \E_{\pi_\theta} \left[ \sum\limits_{t=1}^{H }   \left(\sum\limits_{y\in S} \mathcal{P}( s^{ }_{t,I_t},a^{ }_{t,I_t}, y) V_\xi(y) - V _\xi(  s^{ }_{t,I_t}) + 
   \sum\limits_{i=1}^{I_t-1} \Big(    V_\xi(s_{t, i+1}) - V _\xi(  s_{t,i})\Big) \right)      \right]\\
     &= \E_{\pi_\theta} \left[
 \sum\limits_{t=1}^{H}  
   \sum\limits_{i=1}^{I_t}  c(s_{t,i},a_{t,i})    \right] +\E_{\pi_\theta} \left[ \sum\limits_{t=1}^{H }   \left(\sum\limits_{y\in S} \mathcal{P}( s_{t,I_t},a_{t,I_t}, y) V_\xi(y) - V _\xi(  s_{t,I_t}) + 
   \sum\limits_{i=1}^{I_t-1} \Big(    V_\xi(s_{t, i+1}) - V _\xi(  s_{t,i})\Big) \right)      \right]\\
   &=\E_{\pi_\theta} \left[
 \sum\limits_{t=1}^{H }   \left(c(s_{t,I_t}^{ },a^{ }_{t,I_t}) + \sum\limits_{y\in S}\mathcal{P}( s^{ }_{t,I_t},a^{ }_{t,I_t}, y) V_\xi(y) - V _\xi(  s_{t,I_t}) + 
   \sum\limits_{i=1}^{I_t-1} \Big(c(s_{t,i},a_{t,i}) +   V_\xi(s_{t, i+1}) - V _\xi(  s_{t,i})\Big) \right)      \right]\\
     &= \E_{\pi_\theta} \left[\sum\limits_{t=1}^{H} \sum\limits_{i=1}^{I_t} A_\xi (s_{t,i}, a_{t,i}) \right].
\end{align*}

\end{proof}
\end{lemma}

We define an occupation measure   of policy $\pi_\theta^\Sigma$ at epoch $t$ as a distribution over states  of  $S^\Sigma$:
\begin{align*}
   \mu_\theta(t, s) := \Prob(s_{t} = s),\quad \text{ for each }t=1,..., H,~s\in S^\Sigma,
\end{align*}
 where $s_{t }$ is a state of the  MDP at epoch $t$   under policy $\pi^\Sigma_\theta$.

We define another occupation measure for the states of the SDM process under policy $\pi_\theta$, $\theta\in \Theta$. We denote  the  probability that the SDM process under policy $\pi_\theta$ is at state $y$  after  $i-1$  steps   at epoch $t$, conditioning on $s_{t,1}=s$, as
\begin{align*}
   \phi_\theta(t,i, s, y )  := \Prob(s_{t,i} = y~|~s_{t,1} = s)\quad \text{ for each } t=1,...,H,~i=1,...,I_t,~y\in S,~s\in S^\Sigma. 
\end{align*}

Then, from Lemma \ref{lem:per_diff} we get
\begin{align}\label{eq:two_terms}
      &V_\theta( s_{1, 1}) - V_\xi(s_{1, 1})  \geq   \E_{\pi_\xi} \left[ \sum\limits_{t=1}^{H}\sum\limits_{i=1}^{I_t}\frac{\pi_\theta(a_{t,i}|s_{t,i})}{\pi_\xi(a_{t,i|}s_{t,i})} A_\xi(s_{t,i}, a_{t,i})  \right] \\
      &\quad\quad-\sum\limits_{t=1}^{H }\sum\limits_{s\in  S^\Sigma} \mu_\xi(t, s )  \sum\limits_{i=1}^{I_t}  \max\limits_{s\in S,a\in\A} |A_\xi (s,a)|\sum\limits_{y\in S}|\phi_\xi(t, i, s_t, y) - \phi_\theta(t, i, s_t, y)| \nonumber \\
      &\quad\quad-\sum\limits_{t=1}^{H}\max\limits_{s\in S^\Sigma,~a\in  \A^\Sigma}| A_\xi^\Sigma ( s, a)|\sum\limits_{s\in S^\Sigma} |\mu_\xi(t, s) - \mu_\theta(t,s) |.\nonumber
 \end{align}

Indeed, 
\begin{align*} 
    V_\theta( s_{1,1}) - V_\xi( s_{1, 1})   & =  \E_{\pi_\theta} \left[\sum\limits_{t=1}^{H} \sum\limits_{i=1}^{I_t} A_\xi (s_{t,i}, a_{t,i}) \right]   \\
    & =  \E_{\pi_\theta} \left[\sum\limits_{t=1}^{H}   A_\xi^\Sigma (s_{t }, a_{t}) \right]   \\
     &= \sum\limits_{t=1}^{H }\sum\limits_{s\in  S^\Sigma } \mu_\theta(t, s)  \sum\limits_{a\in \A^\Sigma} \pi^\Sigma_\theta(a| s) A^\Sigma_\xi (s, a)  \\
     & \geq \sum\limits_{t=1}^{H }\sum\limits_{s\in  S^\Sigma} \mu_\xi(t, s)  \sum\limits_{a\in \A^\Sigma} \pi^\Sigma_\theta(a| s) A^\Sigma_\xi (s, a)  \\
     &\quad\quad-\sum\limits_{t=1}^{H}\max\limits_{s\in S^\Sigma,~a\in  \A^\Sigma}| A_\xi^\Sigma ( s, a)|\sum\limits_{s\in S^\Sigma} |\mu_\xi(t, s) - \mu_\theta(t,s) |   \\
     & =   \sum\limits_{t=1}^{H }\sum\limits_{s\in  S^\Sigma} \mu_\xi(t, s)  \sum\limits_{i=1}^{I_t} \sum\limits_{y\in S}\phi_\xi(t, i, s, y) \sum\limits_{a_{t,i}\in \A} \pi_\theta(a_{t,i}| y) A_\xi (y, a_{t,i})  \\
     &\quad\quad-\sum\limits_{t=1}^{H}\max\limits_{s\in S^\Sigma,~a\in  \A^\Sigma}| A_\xi^\Sigma ( s, a)|\sum\limits_{s \in S^\Sigma} |\mu_\xi(t, s) - \mu_\theta(t,s) |   \\
     & \geq \sum\limits_{t=1}^{H }\sum\limits_{s\in  S^\Sigma} \mu_\xi(t, s)  \sum\limits_{i=1}^{I_t} \sum\limits_{y\in S}\phi_\xi(t, i, s, y) \sum\limits_{a_{t,i}\in \A} \pi_\theta(a_{t,i}| y) A_\xi (y, a_{t,i})  \\
     &\quad\quad-\sum\limits_{t=1}^{H }\sum\limits_{s \in  S^\Sigma} \mu_\xi(t, s)  \sum\limits_{i=1}^{I_t}  \max\limits_{s\in S,a\in \A} |A_\xi (s,a)|\sum\limits_{y\in S}|\phi_\xi(t, i, s, y) - \phi_\theta(t, i, s, y)|  \\
     &\quad\quad-\sum\limits_{t=1}^{H}\max\limits_{s\in S^\Sigma,~a\in  \A^\Sigma}| A_\xi^\Sigma ( s, a)|\sum\limits_{s\in S^\Sigma} |\mu_\xi(t, s) - \mu_\theta(t,s) |   \\
     & =\E_{\pi_\xi} \left[ \sum\limits_{t=1}^{H}\sum\limits_{i=1}^{I_t}\frac{\pi_\theta(a_{t,i}|s_{t,i})}{\pi_\xi(a_{t,i|}s_{t,i})} A_\xi(s_{t,i}, a_{t,i})  \right]  \\
     &\quad\quad-\sum\limits_{t=1}^{H }\sum\limits_{s\in  S^\Sigma} \mu_\xi(t, s )  \sum\limits_{i=1}^{I_t}  \max\limits_{s\in S,a\in \A} |A_\xi (s,a)|\sum\limits_{y\in S}|\phi_\xi(t, i, s, y) - \phi_\theta(t, i, s, y)|  \\
     &\quad\quad-\sum\limits_{t=1}^{H}\max\limits_{s\in S^\Sigma,~a\in  \A^\Sigma}| A_\xi^\Sigma ( s, a)|\sum\limits_{s\in S^\Sigma} |\mu_\xi(t, s) - \mu_\theta(t,s) | .
 \end{align*}

Next, we assume that randomized control policy $\pi_\xi$ is the centralized planner's current policy. We want to improve it and get policy $\pi_\theta$  that outperforms the current policy (i.e. $V_\theta(s_{1, 1}) - V_\xi( s_{1,1})>0$).
  We can guarantee the improvement if we find policy $\pi_\theta$ such that the right-hand side (RHS) of   (\ref{eq:two_terms}) is positive. We address the maximization of the RHS of (\ref{eq:two_terms}) following the approach previously used in \cite{Kakade2002, Schulman2015, gluzdai2020}: we  bound $\sum\limits_{s\in \mathcal{S}} |\mu_\xi(t, s) - \mu_\theta(t,s) |$ and $\sum\limits_{y\in S}|\phi_\xi(t,i,s_t,y) - \phi_\theta(t,i,s_t,y)|$   by controlling 
 the maximum change between policies $\pi_\theta$ and $\pi_\xi$, and focus on maximization of the first term of the RHS of (\ref{eq:two_terms}).

In \cite{Schulman2017} the authors proposed  to maximize a \emph{clipping} surrogate objective function:
\begin{align}\label{eq:PO}
L(\theta, \xi):= \E_{{\pi_\xi}} \Big[ \sum\limits_{t=1}^{H} \sum\limits_{i=1}^{I_t}\min& \Big( r_{\theta, \xi}(s_{t,i}, a_{t,i}) A_{\xi} (s_{t,i}, a_{t,i}) , \\
&\text{clip} (r_{\theta, \xi}(s_{t,i}, a_{t,i}),  1-\epsilon, 1+\epsilon)  A_{\xi} (s_{t,i}, a_{t,i})  \Big) \Big],\nonumber 
\end{align}
where $r_{\theta, \xi}(s, a) := \frac{\pi_\theta( s, a)}{\pi_\xi( s, a)}$,   clipping function is equal to
\begin{align*}
  \text{clip}(c,  1-\epsilon, 1+\epsilon):= \begin{cases} 1-\epsilon,~\text{ if } c<1-\epsilon,\\   1+\epsilon,~\text{ if } c>1+\epsilon,\\ c, ~\text{otherwise,}\end{cases}
\end{align*}
and $\epsilon\in(0, 1)$ is a hyperparameter.

We note that the clipping term $\text{clip} (r_{\theta, \xi}( s, a),  1-\epsilon, 1+\epsilon)  A_{\xi}  (  s, a) $ of the objective function  (\ref{eq:PO}) prevents large changes to the policy and keeps $r_{\theta, \xi}(  s, a)$ close to 1, therefore  promoting a conservative update.  
 For the theoretical justification of the PPO algorithm for MDP problems with infinite-horizon discounted reward objectives, see   \cite{Schulman2017},  and for the long-run average cost objectives, see \cite{gluzdai2020}. We note that  \cite{Schulman2017} was motivated by \cite{Schulman2015}, which was inspired by the
conservative policy update ideas pioneered in \cite{Kakade2002}.

We use Monte Carlo simulation to obtain an estimate of the objective function (\ref{eq:PO}). 
We generate $K$ episodes, $K\geq 1$, each  $H$ epochs long. For now, we assume that estimates $\hat A$ of the advantage function required to evaluate  (\ref{eq:PO}) are available. At each step of the SDM process we record a separate datapoint with the following fields (state, action, and advantage function estimate for the state-action pair) to get a dataset:
\begin{align}\label{eq:data}
D_\xi^{(K)}:=  \Big\{  \Big (  s_{t,1,k}, a_{t,1,k}, \hat A (  s_{t,1,k}, a_{t,1,k})  \Big ),\cdots,  
   (   s_{t, I_{t,k} ,k}, a_{t, I_{t,k} ,k}, \hat A ( s_{t, I_{t,k},k}, a_{t,I_{t,k},k})  )  \Big)_{t=1}^{H}
\Big\}_{k=1}^{K}, 
\end{align}
where $s_{t, i, k}$ and $a_{t, i, k}$ are the state  and   action at  the SDM process step $i$, epoch $t$, episode $k$,  respectively.  

Given dataset (\ref{eq:data}) we estimate the objective function as:
 \begin{align}\label{eq:POest}
 \hat L(\theta, \xi, D^{(K)}_\xi):= 
 \frac{1}{K}\sum\limits_{k=1}^K \Big[ \sum\limits_{t=1}^{H} \sum\limits_{i=1}^{I_{t,k}} \min& \Big( r_{\theta, \xi}(s_{t,i,k}, a_{t,i,k}) \hat A_{\xi} (  s_{t,i,k}, a_{t,i,k}) ,   \\
&\quad\text{clip} (r_{\theta, \xi}(  s_{t,i,k}, a_{t,i,k}),  1-\epsilon, 1+\epsilon)  \hat A_{\xi}  (  s_{t,i,k}, a_{t,i,k})  \Big)\Big].\nonumber
\end{align}

Next, we discuss estimating the advantage function of policy $\pi_\xi$. First, we estimate the value function $V_\xi$. We compute a Monte Carlo estimate of the value function that corresponds to each step in the generated episodes (\ref{eq:data}), such as
\begin{align}\label{eq:Vest}
    \hat V_{t, i, k} := \sum\limits_{j=i}^{I_{t,  k}} c(s_{t,j,k}, a_{t,j,k}) +\sum\limits_{\ell = t+1}^{H}  \sum\limits_{j=1}^{I_{\ell,k}} c(s_{\ell,j, k}, a_{\ell,j, k}),
\end{align}
which is a one-replication estimate of the value function $V(s_{t, i, k})$ at state $s_{t, i, k}$ that is visited at epoch $t$, episode $k$, after  $i-1$ steps of the SDM process.

We   use function approximator $V_\psi:S\rightarrow \R$ to get a low-dimentional representation of value function $V_\xi$. We consider a set of function approximators $\{V_\psi, \psi\in \Psi\}$, and based on one-replication estimates we find the optimal $V_\psi $ that minimizes the mean-square norm:
 \begin{align}\label{eq:mse}
     \sum\limits_{k=1}^{K}\sum\limits_{t=1}^{H}\sum\limits_{i=1}^{I_{t,k}} \| V_\psi (s_{t, i, k}) - \hat V_{t, i, k} \|^2.
 \end{align}

Next, we obtain the advantage function estimates
\begin{align}\label{eq:Aest}
    \hat A ( s_{t, i,  k}, a_{t, i, k}):= 
     \begin{cases}
    c(s_{t, i, k}, a_{t, i, k}) +V_\psi( s_{t, i+1, k}) - V_\psi(s_{t, i, k})\quad \text{ if } i\neq I_{t,  k},\\
    c (s_{t, i, k}, a_{t, i, k}) +V_\psi (s_{t+1, 1, k}) - V_\psi(s_{t, i, k})\quad \text{ otherwise}, 
\end{cases}
\end{align}
for each $t=1, ...,H$; $k=1,..., K;$ and $i=1,..., I_{t,k}$.

 Our proposed PPO algorithm consists of the following steps.
 
\begin{algorithm}[H]
\SetAlgoLined
\LinesNumbered
\SetKwBlock{Begin}{Begin}{}
\SetAlgoLined
\SetKwProg{Loop}{LOOP}{}{}
\KwResult{policy $\pi_{\theta_J}$} 
 Initialize  policy function  $\pi_{\theta_0}$ and value function approximator $V_{\psi_{-1}}$\;
 \For{ policy iteration $j = 1, 2, \dotsc, J$}{
  
Run policy $\pi_{\theta_{j-1}}$ for $K$ episodes and collect dataset (\ref{eq:data}).
 
   Construct  Monte-Carlo estimates of the value function $V_{\theta_{j-1}}$  following (\ref{eq:Vest}).

Update function approximator $V_\psi$ minimizing (\ref{eq:mse}).  
        
Estimate advantage functions $\hat A (s_{t, i,  k}, a_{t, i, k})$ by (\ref{eq:Aest}).
      
 Maximize surrogate objective function (\ref{eq:POest}) w.r.t. $\theta$.
 Update $\theta_{j} \leftarrow  \theta$
      }
 \caption{The PPO algorithm} 
 \label{alg:ppo}
\end{algorithm}

 \section{Experiment and Results}\label{sec:num-study}

In this section we evaluate the performance of our proposed PPO for a ride-hailing service's transportation network consisting of five regions, $N = 1000$ cars, and $H = 360$ minutes, based on real data from Didi Chuxing.

 Following \cite{BravDaiLiuLei2019}, at the start of each working day, the centralized planner distributes the cars in proportion to each region's expected demand.
 We set patience time at $L = 5$. We establish the reward functions (i.e., car-passenger matching rewards are equal to $c^f_t(o,d,\eta) = 1$, and empty-car routing costs are equal to $c_t^e(o,d) = 0,$ for each $o, d = 1,\dotsc,R$, $\eta = 0,1,\dotsc,L $, and $t = 1,\dotsc, H$) such that the total reward accumulated at the end of the working day corresponds to the number of completed ride requests. The remaining experiment details can be found in Appendix~EC.3.2. of ~\cite{BravDaiLiuLei2019}.




In this way, the total reward accumulated by the end of a working day correspond to the  number  of ride requests fulfilled. This can be reinterpreted as the fraction of ride requests fulfilled,  given  a sample path of the passenger arrivals. The number of completed ride requests fulfilled is the common objective considered by the dynamic matching problems, see, for example \cite{OzkaWard2017}.

We use two separate and fully connected feed-forward neural networks to represent randomized control policies $\pi_\theta$,  $\theta\in \Theta$ and value functions $V_\psi$, $\psi\in \Psi$. We run the algorithm for $J =75$ policy iterations. The algorithm simulates $K = 300$ episodes in each iteration.  See Appendix for the details.
We use the best result for the five-region experiment in \cite{BravDaiLiuLei2019} as the benchmark. The ``time-dependent lookahead'' policy from \cite{BravDaiLiuLei2019} could achieve 84\% fulfilled ride requests.

\begin{figure}[thpb]
      \centering
      \includegraphics[scale=0.7]{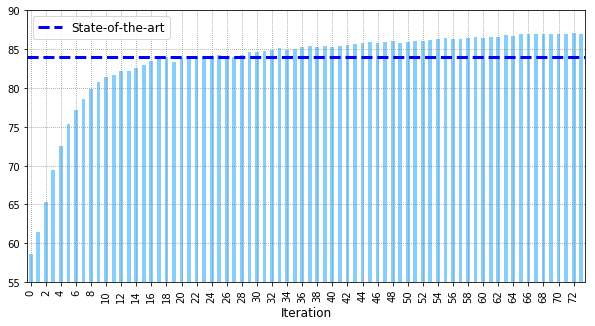}
      \caption{Learning curve from Algorithm \ref{alg:ppo} for the   transportation network from \cite{BravDaiLiuLei2019}. The dashed blue line shows the
best performance of the ``time-dependent lookahead'' policy from \cite{BravDaiLiuLei2019}. The blue columns  show the  performance of the randomized control policies obtained at every iteration of Algorithm \ref{alg:ppo}.  }
  \label{fig:average_matching_rates}
  \end{figure}

Figure \ref{fig:average_matching_rates} shows that our PPO algorithm   achieves  $87\%$  fulfilled ride requests  after $J = 75$ policy iterations. In fact,  it needs only $8$ policy iterations to boost the performance to $80\%$ from the initial $59\%$ using a policy NN with random weights. To  evaluate the performance of the randomized control policy after every iteration,  we run the policy for $K=300$ working days, and take the average of the fractions of fulfilled ride requests on each day.

\section{Conclusion}
\label{sec:conclusion}
This paper proposes a method to optimize the  finite horizon total reward objective in a ride-hailing service system.
The large  action space 
prohibits the direct use of policy optimization RL methods. Although the standard PPO algorithm \cite{Schulman2017} suggests designing a policy NN such that the the number of units in its output layer is equal to the number of actions, the algorithm becomes computationally expensive because the parameters grow exponentially with the number of agents, and it cannot capture the similarities among the actions because it has to consider them separately.

Instead of searching for the optimal actions directly, we use the PPO algorithm to learn the most beneficial type of trip
to fulfill at a given system state. Repeated executions of the PPO policy allow a centralized planner to prioritize 
trips and sequentially assign tasks to all available cars.

  A numerical experiment using real data from a ride-hailing service demonstrates that the PPO algorithm applied to our MDP formulation outperforms  the
policy proposed in \cite{BravDaiLiuLei2019} by 3\%.

\section*{Acknowledgment}  This research  was  supported in part by   National Science Foundation Grant CMMI-1537795.

\bibliographystyle{plain}
\bibliography{ride-hailing_paper}

\begin{thebibliography}{10}

\bibitem{BravDaiLiuLei2019}
A.~Braverman, J.~G. Dai, X.~Liu, and Y.~Lei.
\newblock {Empty-car routing in ridesharing systems}.
\newblock {\em Operations Research}, 67(5):1437--1452, 2019.

\bibitem{gluzdai2020}
J.~G. Dai and M.~Gluzman.
\newblock Queueing network controls via deep reinforcement learning.
\newblock 2020.
\newblock \url{https://arxiv.org/abs/2008.01644}.

\bibitem{Henderson2002}
Shane~G. Henderson and Peter~W. Glynn.
\newblock {Approximating martingales for variance reduction in Markov process
  simulation}.
\newblock {\em Mathematics of Operations Research}, 27(2):253--271, 2002.

\bibitem{Kakade2002}
S.~Kakade and J.~Langford.
\newblock {Approximately Optimal Approximate Reinforcement Learning}.
\newblock In {\em Proc. ICML'02}, pages 267--274, 2002.

\bibitem{KeXiao19}
J.~Ke, F.~Xiao, H.~Yang, and J.~Ye.
\newblock {Optimizing Online Matching for Ride-Sourcing Services with
  Multi-Agent Deep Reinforcement Learning}.
\newblock 2019.
\newblock \url{https://arxiv.org/abs/1902.06228}.

\bibitem{Mnih2015}
V.~Mnih, K.~Kavukcuoglu, D.~Silver, A.~A. Rusu, and et. al.
\newblock {Human-level control through deep reinforcement learning}.
\newblock {\em Nature}, 518(7540):529--533, 2015.

\bibitem{Oda2018}
T.~Oda and C.~Joe-Wong.
\newblock Movi: A model-free approach to dynamic fleet management.
\newblock In {\em IEEE INFOCOM'18}, pages 2708--2716, 2018.

\bibitem{OzkaWard2017}
E.~Ozkan and A.~Ward.
\newblock Dynamic matching for real-time ride sharing.
\newblock {\em Stochastic Systems}, 10(1):29--70, 2020.

\bibitem{Schlobach2018}
M.~Schlobach and S.~Retzer.
\newblock {Didi Chuxing - How China's ride-hailing leader aims to transform the
  future of mobility}, 2018.

\bibitem{Schulman2015}
J.~Schulman, S.~Levine, P.~Moritz, M.~I. Jordan, and P.~Abbeel.
\newblock {Trust Region Policy Optimization}.
\newblock In {\em ICML'15}, pages 1889--1897, 2015.

\bibitem{Schulman2017}
J.~Schulman, F.~Wolski, P.~Dhariwal, A.~Radford, and O.~Klimov.
\newblock {Proximal Policy Optimization Algorithms}.
\newblock 2017.
\newblock \url{http://arxiv.org/abs/1707.06347}.

\bibitem{ShiGao19}
J.~Shi, Y.~Gao, W.~Wang, N.~Yu, and P.~A. Ioannou.
\newblock {Operating Electric Vehicle Fleet for Ride-Hailing Services With
  Reinforcement Learning}.
\newblock {\em IEEE Transactions on Intelligent Transportation Systems}, pages
  1--13, 2019.

\bibitem{Tang2019ADV}
X.~Tang, Z.~Qin, F.~Zhang, Z.~Wang, Z.~Xu, Y.~Ma, H.~Zhu, and J.~Ye.
\newblock A deep value-network based approach for multi-driver order
  dispatching.
\newblock In {\em KDD '19}, page 1780–1790, 2019.

\bibitem{WangQinTanYeZhu18}
Z.~Wang, Z.~Qin, X.~Tang, J.~Ye, and H.~Zhu.
\newblock {Deep Reinforcement Learning with Knowledge Transfer for Online Rides
  Order Dispatching}.
\newblock In {\em ICDM'18}, pages 617--626, 2018.

\bibitem{Xu2018ADV}
Z.~Xu, Z.~Li, Q.~Guan, D.~Zhang, Q.~Li, J.~Nan, C.~Liu, W.~Bian, and J.~Ye.
\newblock Large-scale order dispatch in on-demand ride-sharing platforms: A
  learning and planning approach.
\newblock In {\em KDD'18}, page 905–913, 2018.

\end{thebibliography}

  \section*{Appendix (details of the experiment setting in Section \ref{sec:num-study})}



 We  refer to the neural network used to represent a policy as \textit{the  policy NN} and to  the neural network used to approximate a value function as \textit{the value function NN}.  
Both policy and value NNs consist of an input layer, an embedding layer, three hidden layers, and an output layer. The embedding layer deals with the time-of-day component in a state, which is a categorical variable taking one of $360$ values.
The sizes of the input layer, embedding layer, first -- third hidden layers are the same for both NNs and equal to 394, 6, 399, 44, 5, respectively. 
The output layer of the value NN has one unit using identity activation function. The output layer of the policy NN has 25 units using softmax activation function.

 Table \ref{tab:dNN_hyperparams} summarizes the hyper-parameters of Algorithm \ref{alg:ppo} used in the 5-region experiment in Section \ref{sec:num-study}. For the policy NN learning rate and clipping, we adopt a simple decay scheme as follows: 
\begin{align*}
    \beta_j  := \max(1-j/J, 0.01) \beta ,\quad
    \epsilon_j  := \max((1-j/J) \epsilon, 0.01),
\end{align*}
where $\beta,$ $\epsilon$ denote the initial learning rate and clipping reported in Table \ref{tab:dNN_hyperparams}, and $j$ denotes the $j$-th policy iteration.

\begin{table}[h!]
    \centering
    \begin{tabular}{|c|c|c|}
    \hline 
         Hyper-parameter & Description & Value \\
         \hline
         $J$& Number of policy iterations & $75$ \\
         \hline 
         $K$&Number of episodes per policy iteration & $300$ \\
         \hline 
         $\beta $&Initial learning rate for the policy NN &
         $0.00005$\\ 
         \hline 
         $-$&Learning rate for the value NN &
         $0.0001$ \\
         \hline 
         $\epsilon$&Initial clipping & $0.2$ \\
         \hline 
         $-$ & Number of passes over training data for policy NN update & $3$\\
         \hline 
         $-$ & Number of passes over training data for value NN update & $10$\\
         \hline 
         $-$ & \begin{tabular}{@{}c@{}}Kullback–Leibler (KL) target for prompting early stopping  \\ during policy NN training
         \end{tabular}
         & $0.012$\\
         \hline 
         $-$ & $L_2$ regularization factor for the embedding layers & $0.005$\\
         \hline 
    \end{tabular}
    \caption{Hyper-parameters for Algorithm \ref{alg:ppo}.} 
    \label{tab:dNN_hyperparams}
\end{table}

Table \ref{tab:mdp-params} reports the values for the parameters $ \lambda,$ $P$, and $\tau$ that characterize the time-varying traffic patterns of the five-region transportation network considered in Section~\ref{sec:num-study}. These values (except the patient time $L$) are the same as in Appendix~EC.3.2. of ~\cite{BravDaiLiuLei2019}.

\begin{table}[h!]
    \centering
    \begin{tabular}{|c|c|c|}
    \hline 
         Parameter & Description & Value \\
         \hline
         $R$ & Number of regions & $5$ \\
         \hline 
         $N$ & Number of cars & $1000$ \\ 
         \hline 
         $H$ & Length of a working day (minute) &$360$ \\
         \hline 
         $L$ & Passenger patience time (minute) & $5$ \\
         \hline 
         $ \lambda $ &
         Passenger arrival rates at each region per minute 
         &   \eqref{eq:traf-pat-1} --~\eqref{eq:traf-pat-3} \\
         \hline 
         $ P $ & Probabilities of going to each destination from each origin &  \eqref{eq:traf-pat-1} --~\eqref{eq:traf-pat-3} \\
         \hline 
         $ \tau $ &Mean travel times (minute) &  \eqref{eq:traf-pat-1} --~\eqref{eq:traf-pat-3}\\
         \hline 
         $ c^f $, $c^e$ & Immediate rewards for an atomic action& 
         $c^f(o,d,\eta)\equiv 1$, $c^e(o,d)\equiv 0$.
         \\
         \hline 
         \end{tabular}
         \caption{Traffic parameters for the transportation network.}
         \label{tab:mdp-params}
\end{table}


For the first two hours, $t = 1,2,\dotsc, 120$, traffic parameters are equal to
\begin{align}
     \lambda(t) 
    &=
    \begin{bmatrix}
    1.8\\
    1.8\\
    1.8\\
    1.8\\
    18
    \end{bmatrix}
    ,\quad 
     P(t)
    = \begin{bmatrix}
    0.6 & 0.1 & 0 & 0.3 & 0\\
    0.1 & 0.6 &0  &0.3  &0\\
    0& 0 &0.7& 0.3 &0 \\
    0.2 &0.2& 0.2& 0.2 &0.2 \\
    0.3& 0.3& 0.3& 0.1& 0
    \end{bmatrix}
    , \quad 
     \tau(t)  = \begin{bmatrix}
    9 &15 &75 &12 &24\\
    15  &6 &66 & 6 &18 \\
    75 &66 & 6 &60 &39 \\ 
 15 & 9 &60  &9 &15 \\ 
 30 &24 &45 &15 &12
    \end{bmatrix}, 
    \label{eq:traf-pat-1}
\end{align}
where  $\lambda_o(t)$ is the expected  number of passengers arriving to region $o$ at time $t$,  $ P_{od}(t)$  is a probability that  a passenger arriving at region $o$ at time $t$ requests a trip to destination region $d$, and $\tau_{od}(t)$ is a duration of the trip, started at time $t$, from region $o$ to region $d$.

For the third, fourth hours, $t = 121, 122, \dotsc, 240$, the traffic parameters are equal to
\begin{align}
     \lambda(t) 
    &=
    \begin{bmatrix}
    12\\
    8\\
    8\\
    8\\
    2
    \end{bmatrix},
    \quad 
     P(t)
    = \begin{bmatrix}
   0.1 &0 &0 &0.9 &0 \\
   0& 0.1& 0 &0.9 &0 \\
   0 &0 &0.1& 0.9& 0 \\
   0.05 & 0.05 & 0.05 & 0.8 & 0.05 \\
   0 &0 &0& 0.9 &0.1 
    \end{bmatrix}
        , \quad 
     \tau(t)  = \begin{bmatrix}
    9 &15& 75& 12& 24\\
 15 & 6& 66&  6 &18\\
 75 &66 & 6& 60 &39\\
 12&  6& 60 & 9 &15\\
 24 &18 &39& 15& 12\\
    \end{bmatrix}.
    \label{eq:traf-pat-2}
\end{align}

For the last two hours, $t = 241, 242, \dotsc, 360$,  the traffic parameters are equal to
\begin{align}
     \lambda(t) 
    &=
    \begin{bmatrix}
    2\\
    2\\
    2\\
    22\\
    2
    \end{bmatrix},\quad 
     P(t)
    = \begin{bmatrix}
   0.9 &0.05 &0 &0.05 &0\\
   0.05 &0.9 &0 &0.05 &0 \\
   0 &0 &0.9 &0.1 &0 \\ 
   0.3 &0.3 &0.3 &0.05 &0.05 \\
   0 &0 &0 &0.1& 0.9  
    \end{bmatrix} , \quad 
     \tau(t)  = \begin{bmatrix}
   9& 15& 75& 12& 24\\
 15&  6& 66&  6& 18\\
 75& 66&  6& 60& 39\\
 12&  6& 60&  9& 15\\
 24& 18& 39& 15& 12\\
    \end{bmatrix}.
    \label{eq:traf-pat-3}
\end{align}
\end{document}